\documentclass[a4paper,twoside]{article}
\usepackage{amsfonts,amssymb,mathrsfs,amsmath}
\title{POLYBOXES, CUBE TILINGS AND RIGIDITY}
\date{}
\author{Andrzej P. Kisielewicz and Krzysztof Przes{\l}awski\\
\\
{\small Wydzia{\l} Matematyki, Informatyki i Ekonometrii, Uniwersytet Zielonog\'orski}\\
{\small ul. Z. Szafrana 4a, 65-516 Zielona G\'ora, Poland}\\
{\small A.Kisielewicz@wmie.uz.zgora.pl}\\
{\small K.Przeslawski@wmie.uz.zgora.pl}}
\newtheorem{pr}{\sc Proposition}
\newtheorem{lemat}[pr]{\sc Lemma}
\newtheorem{tw}[pr]{\sc Theorem}
\newtheorem{wn}[pr]{\sc Corollary}
\newtheorem{df}{\sc Definition}

\newtheorem{uw}{\rm REMARK}
\newtheorem{uwi}[uw]{\rm REMARKS}
\newenvironment{uwa}{\begin{uw} \rm}{\end{uw}}

\newtheorem{nap}{\rm EXAMPLE}

\newtheorem{nps}[nap]{\rm EXAMPLES}

\def\ka #1{\mathscr{#1}}
\def\kal #1 #2{\mathscr{#1}^{#2}}
\def\proof{\noindent \textit{Proof.\,\,\,}}
\def\dowod #1{\noindent\textit{Proof #1.\,\,\,}}

\def\add #1{\mathbf{A}(#1)}
\def\addb #1{\mathbf{A}_{\mathrm{B}} (#1)}
\def\modgen #1{\zet \pud #1}
\def\Modgen #1{\zet \Pud #1}

\def\zet{\mathbb{Z}}

\def\er{\mathbb{R}}

\def\te{\mathbb{T}}
\def\pud #1{\operatorname{box}(#1)}
\def\Pud #1{\operatorname{Box}(#1)}
\def\pudd #1{\operatorname{box}^2(#1)}

\def\iver #1{\mbox{\tt [} #1 \mbox{\tt]}}
\def\index #1#2{\mathrm{ind}\,(#1,#2)}
\def\Index #1#2{\mathrm{Ind}\,(#1,#2)}
\def\supp{\operatorname{supp}}
\def\przyst #1{(\operatorname{mod} #1)}
\def\naw #1{^{(#1)}}
\begin{document}
\maketitle
\begin{abstract}
A non-empty subset $A$ of $X=X_1\times\cdots \times X_d$ is a (proper) box if $A=A_1\times \cdots\times A_d$ and $A_i\subset X_i$ for each $i$. Suppose that for each pair of boxes $A$, $B$ and each $i$, one can only know 
which of the three states takes place: $A_i=B_i$, $A_i=X_i\setminus B_i$, $A_i\not\in\{B_i, X_i\setminus B_i\}$. Let $\ka F$ and $\ka G$ be two systems of disjoint boxes. Can one decide whether $\bigcup \ka F=\bigcup \ka G$? In general, the answer is `no', but as is shown in the paper, it is `yes' if both systems consist of pairwise dichotomous boxes.
Several criteria that enable to compare such  systems are collected. The paper includes also rigidity results, which 
say what assumptions have to be imposed on $\ka F$ to ensure that $\bigcup \ka F=\bigcup \ka G$ implies $\ka F=\ka G$.
As an application, the rigidity conjecture for $2$-extremal cube tilings of Lagarias and Shor is verified.        


\end{abstract}

\section{Introduction}
Let $\te^d=\{(x_1,\ldots, x_d)(\operatorname{mod}2):(x_1,\ldots ,x_d)\in \er^d\}$ be a flat torus. Suppose that $[0,1)^d+\lambda$, $\lambda\in \Lambda$, is a cube tiling of $\te^d$.  
The cube tiling $[0,1)^d+\Lambda$ is 2-\textit{extremal} if for each $\lambda\in \Lambda$ there is a unique $\lambda'\in \Lambda$ such that $ (|\lambda_1-\lambda'_1|,\ldots , |\lambda_d-\lambda'_d|)\in \zet^d$. Let $\Lambda_+$, $\Lambda_-$ be any decomposition of $\Lambda$ such that each of the component does not contain any of the pairs $\{\lambda, \lambda'\}$, $\lambda\in \Lambda$. In an important paper \cite{LS2}, where cube tilings of $\er^d$ contradicting  Keller's celebrated conjecture (see \cite{K1,K2,P,LS1,Ma,SSz}) in a certain strong sense are constructed, Lagarias and Shor conjectured that $\Lambda_+$ and $\Lambda_-$ determine each other; that is, if $[0,1)^d+\Gamma$ is another 2-extremal cube tiling of $\te^d$, and $\Gamma_+$, $\Gamma_-$ is a corresponding decomposition of $\Gamma$, then the equality $\Lambda_+=\Gamma_+$ implies  $\Lambda_-=\Gamma_-$. (Actually, their assertion, called in \cite{LS2} \textit{the rigidity conjecture for 2-extremal cube-tilings}, is stated in the language of 2$\zet^d$-periodic cube tilings of $\er^d$.) In this paper we show that a far reaching generalization of the rigidity conjecture remains valid (Theorem \ref{krewmleko}). In a sense, we could say that the latter result is a by-product of the present investigations. We arrive at this problem working with slightly different structures: partitions of the Cartesian products of the finite sets, called here boxes, into boxes. Our interest in these structures comes from a certain minimization problem of Kearnes and Kiss \cite{KK}, which has been solved by Alon, Bohman, Holzman and Kleitman \cite{ABHK}.  Minimal partitions that are involved in their solution have been characterized in \cite{GKP}. In Section \ref{minimal} we extend these investigations to what we call polyboxes. 

A non-empty subset $A$ of the Cartesian product
$X:=X_1\times\cdots\times X_d$ of finite sets $X_i$, $i\in[d]$,
is called a \textit{ box} if $A=A_1\times\cdots \times A_d$ and
$A_i\subseteq X_i$ for each $i\in [d]$. We say that $A$ is a $k$-dimensional box if
$| \{i\in [d]: |A_i|>1\}|=k$. We call $A$  \textit{ proper}
if $A_i\neq X_i$ for each $i\in [d]$.  The family of all boxes contained in $X$
is denoted by $\Pud X$, while $\pud X$ stands for the family of all proper boxes
in $X$.

Two boxes $A$ and $B$ in $X$ are said to be \textit{dichotomous}
if there is an $i\in [d]$ such that $A_i=X_i\setminus B_i$. Any collection of pairwise
dichotomous boxes is called a \textit{suit}. A suit is \textit{proper} if consists
of proper boxes. A non-empty set
$F\subseteq X$ is said to be a \textit{ polybox} if
there is a suit $\ka F$ for $F$, that is,  $\bigcup \ka F=F$. 

Since polyboxes 
are defined by means of partitions into boxes, it is not surprising that  characteristics of polyboxes will be expressed in terms of partitions as well. Therefore, such characteristics should be invariant on the choice of a partition.
In Section \ref{additive functions}, we define the class of additive functions on $\pud X$. All characteristics of polyboxes that appear in the paper are defined with the use of additive functions. Theorem \ref{add1} plays in this respect a crucial role.  

An interesting class of invariants is described in Sections \ref{dilabel} and \ref{binary codes}.   

In Section \ref{indices} we discuss important numerical characteristics of a polybox, the indices. 

In the following section we show, among other things, that certain mild assumptions on the symmetry of a polybox 
imply that all  indices of the polybox are even numbers.     

In Section \ref{rikitiki} we give a sufficient condition which guarantees a polybox to be rigid in the sens that it has a unique proper suit (Theorem \ref{sztywnypal}). This result appears again in a greater generality, applicable to the already mentioned case of the rigidity of cube tilings, in Section \ref{words} (Theorems \ref{niejednoznacznosc} and \ref{rikitikigenom}).

One of the  basic questions is whether two given suits define the same polybox. We address this question in several places. An important procedure, which enables us to answer it, is described in Section \ref{modules} (Remark \ref{kryt2}).
This procedure is based on  a certain decomposition of  the free $\zet$-module  generated by boxes (Theorem \ref{rzut}). 
       
The results of the paper are summarized in an abstract setting of words in Section \ref{words}. There is also defined and investigated an interesting cover relation.

\section{Minimal partitions}
\label{minimal}
Let  $F$ be a subset of a  $d$-dimensional box $X$. A partition of $F$ into
proper boxes is \textit{minimal} if it is of minimal cardinality among all such
partitions. It is observed in \cite{GKP} that if $F=X$,
then the minimal partitions of $F$ coincide with the proper suits for $F$.
This result extends to polyboxes:

\begin{tw}
\label{minor}
If $F$ is  a polybox in a $d$-dimensional box $X$ and $\ka F\subseteq \pud X$
is a partition of $F$, then $\ka F$ is minimal if and only if $\ka F$ is a suit.
\end{tw}

The proof is a refinement of an argument given in \cite{ABHK},
and is much the same as in \cite{GKP}, however, we added 
to it a geometric flavour.

Before going into the  proof, we collect several indispensable definitions and lemmas.

Let ${\ka O}(X_i)$ be the family of all sets of odd size which are contained in $X_i$. Let $B$ be a
 subset of $X$. We define $\widehat {B}$ to be the subset of
${\ka O}(X_1)\times \cdots \times{\ka O}(X_d)$ that consists of all
$d$-tuples $(A_1,\ldots , A_d)$ for which the set
$B\cap (A_1\times\cdots \times A_d)$ is of odd size.

Suppose that $B$ is a box. Let $\ka O B_i$ be the set of all sets of odd size that are contained
in $X_i$ such that
their intersections with $B_i$ are of odd size as well.  One can easily observe that $\widehat{B}=
\ka O B_1\times\cdots \times \ka O B_d$.
In particular, $\widehat{B}$ is a box. Moreover,
since
\begin{equation}
\label{os}
|\ka O {B_i}|=\frac{|\ka O (X_i)|}{2}= \frac{2^{|X_i|}}{4}
\end{equation}
for each $i\in [d]$, we obtain
\begin{equation}
\label{kostkan}
|\widehat{B}|=2^{|X|_1-2d},
\end{equation}
where $|X|_1$ is defined by the equation $|X|_1=|X_1|+\cdots +|X_d|$.

\begin{lemat}
\label{comp} The following conditions are equivalent:
\begin{description}
\item[(i)]
boxes $B, C \in \Pud X$ are dichotomous,
\item[(ii)]
$\widehat{B}$ and $\widehat{C}$ are dichotomous,
\item[(iii)]
$\widehat{B}$ and $\widehat{C}$ are disjoint.
\end{description}
\end{lemat}

\proof The equivalence `(i)$\Leftrightarrow $(ii)'
is deduced easily from the observation that $X_i\setminus B_i=C_i$
if and only if $\ka O (X_i)\setminus \ka O B_i=\ka O C_i$.
Concerning  `(ii)$\Leftrightarrow $(iii)', only the implication
`(ii)$\Leftarrow $(iii)' is non-trivial.
Thus, if $\widehat{B}$ and $\widehat{C}$ are disjoint,
then there is an  $i\in [d]$ such that $\ka OB_i$ and $\ka OC_i$ are disjoint.
By (\ref{os}), each of these sets contains
half of the elements of  $\ka O(X_i)$. Therefore, they are complementary.
Consequently, $\widehat{B}$ and $\widehat{C}$ are dichotomous.\hfill $\square$

\medskip
\medskip
The next lemma is rather obvious.
\begin{lemat}
\label{symmdiff}
If $B$ and $C$ are  subsets of  $X$, then
$$
\widehat{B\ominus C}=\widehat{B}\ominus\widehat{C},
$$
where $\ominus $ denotes the symmetric difference.
\end{lemat}

\textit{ Proof of the theorem.\,\,} Let us enumerate all elements of $\ka F$, that is,
$\ka F =\{B^1,\ldots ,B^m\}$. Let $\widehat{\ka F}=\{\widehat{B^1}, \ldots ,
\widehat{B^m}\}$. By the preceding lemma, we have
$$
\widehat{F}=\widehat{B^1}\ominus\cdots\ominus \widehat{B^m}.
$$
This equality implies
$$
|\widehat{F}|\leq |\widehat{B^1}|+\cdots + |\widehat{B^m}|,
$$
where the equality holds if and only if  the elements of $\widehat{\ka F}$
are mutually disjoint.
If we divide the above inequality by $2^{|X|_1-2d}$, then, by (\ref{kostkan}),
we obtain
\begin{equation}
\label{nier1}
\frac{|\widehat{F}|}{2^{|X|_1-2d}}\leq m,
\end{equation}
which means that the size of any partition of $F$ into proper boxes is bounded from
below by the left side of the above inequality. Moreover, this bound is
tight and, according to Lemma \ref{comp}, is attained if and only if
$\ka F$ is a suit.\hfill $\square$

\medskip
\medskip
Observe that, as a by-product, we have shown that the  minimal partitions  of a
polybox $F$ have their size equal to the number standing on the left side of
(\ref{nier1}). This suggests
the following definition: Let $G\subseteq X$. The number $|G|_0$, given by the formula
\begin{equation}
\label{vol0}
|G|_0={|\widehat{G}|}{2^{|X|_1-2d}}=2^d\frac{|\widehat{G}|}{|\widehat{X}|},
\end{equation}
is called the \textit{box number} of $G$.

By much the same method as applied  above, one can obtain the following
characterization of polyboxes.
\begin{tw}
\label{vol00}
Let $G$ be a non-empty subset of a $d$-box $X:=X_1\times \cdots \times X_d$.
The size of a minimal partition of $G$ into proper boxes is at least $|G|_0$. Moreover,
$G$ is a polybox if and only if there is a partition of $G$ into proper
boxes of size $|G|_0$.
\end{tw}
By the definition of the box number,
\begin{equation}
|G|_0 \le |X|_0= 2^d.
\end{equation}
Moreover, equality occurs if and only if $G=X$. Indeed, for $x\in X$, let us define
 $\hat{x}=(\{x_1\},\ldots, \{x_d\})\in \widehat{X}$.
It is easily seen that $x\in X\setminus G$ if and only if
$\hat{x}\in\widehat{X}\setminus\widehat{G}$.
Thus, we have
\begin{tw}
\label{krakow}
Let $F$ be a polybox contained in a $d$-box $X$. If there is a proper suit $\ka F$
for $F$ of size $2^d$, then
$F=X$.
\end{tw}

The following interesting question arises:

\medskip
\noindent {\sc Question.} Suppose that $F$ and $G$ are polyboxes contained in a $d$-box $X$.
Is it true that if $F\subseteq G$, then $|F|_0\le |G|_0?$
\medskip
\begin{pr}
\label{suma}
If $F$ and $G$ are disjoint polyboxes in $X$ and there
are proper suits $\ka F$ for  $F$, and
$\ka G$ for $G$  such that $\ka F\cup\ka G$ is a suit for $ F\cup G$,
then the same holds true for every pair of proper suits for $F$ and $G$.
\end{pr}

\proof Let $\ka F^*$ and $\ka G^*$ be proper suits for $F$ and $G$, respectively.
By Theorem~\ref{minor}, proper suits for the same polybox are of the same
size. Consequently,
$$
|\ka F^*\cup \ka G^*|=|\ka F^*|+|\ka G^*|=|\ka F|+|\ka G|= |\ka F\cup\ka G|.
$$
Again by Theorem \ref{minor}, the fact that $\ka F\cup \ka G$ is a suit, and the
preceding equality, we conclude that the partition
$\ka F^*\cup\ka G^*$ is a suit.\hfill$\square$

If sets $F$ and $G$ are as described in Proposition \ref{suma}, then we call them \textit{strongly disjoint}. 

\section{Additive functions}
\label{additive functions}
Let a box $A=A_1\times \cdots \times A_d$ and a set
$I\subseteq [d]$ be given. Let  $i_1,\ldots, i_k$ be the elements of $I$
written in increasing order. We define
$A_I:=A_{i_1}\times\cdots \times A_{i_k}$. We have the natural projection
$a\mapsto a_I$ from
$A$ onto $A_I$, where if $a=(a_1,\ldots , a_d)$, then $a_I=(a_{i_1},\ldots , a_{i_k})$.
If $B\subseteq A$, then we put $B_I=\{b_I:b\in B\}$, and if $\ka F\subseteq 2^A$, then $\ka F_I=\{B_I: B\in\ka F\}$. 
To simplify our notation, we shall write $i'$ rather than $[d]\setminus \{i\}$. 

We say that two boxes $A$ and $B$ contained in a $d$-box $X$ form a
\textit{twin pair} 
if
$A_{i'}=B_{i'}$ and $A_i=X_i\setminus B_i$.

If $p(x)$ is a sententional function, then, as proposed by Iverson, $\iver {p(a)}=1$, if $p(a)$ is true for $x=a$,
and $\iver{p(a)}=0$, if $p(a)$ is false for $x=a$.

Let us extend the notation introduced in Section \ref{minimal} letting $\ka O(X)\subset \Pud X$
be the family of all boxes of odd size contained in $\Pud X$.

Let $X$ be a $d$-box and let $M$ be a module over a commutative ring $R$.
A function $f\colon \pud X\to M$ is  \textit{additive} if for any two twin pairs $A$, $B$ and $C$, $D$,
the equation $A\cup B=C\cup D$ implies
$$
f(A)+f(B)=f(C)+f(D).
$$
The module of all $M$-valued additive  functions is denoted $\add{X,M}$ .

In this paper we shall be concerned with real valued additive functions.

For each $B\in\ka O (X)$,  let us define $\eta_B\colon 2^{X}\to \er$ by
$$
\eta_B(A)=\iver{ |A\cap B|\equiv 1  \mod{2}}.
$$
It is clear that the restriction of $\eta_B$ to $\pud X$ is additive.

Now, let us confine ourselves to the case $d=1$, that is, we shall assume that $X$
is simply a finite set that contains at least two elements. Then $\pud X$ coincides with
$2^X\setminus\{X, \emptyset\}$.

For each $C\in \pud X$, let $\varphi_C \colon\pud X\to \er$ be defined by
\begin{equation}
\label{fione}
\varphi_C (A)=\iver{A=C}-\iver{A=X\setminus C}.
\end{equation}
Moreover, let $\varphi_X \colon \pud X\to  \er$ be the constant function equal to 1.
It is clear that the  functions $\varphi_C$, $C\in \Pud X$, are additive.
It is also clear that $\varphi_C\perp\varphi_D$,
whenever $C\not\in\{D, X\setminus D\}$, where the orthogonality is related to the scalar product defined by
\begin{equation}
\label{skalar}
\langle f, g\rangle = \sum_{C\in \pud X} f(C)g(C).
\end{equation}
Let us note for future reference that
\begin{equation}
\label{sym}
\varphi_{X\setminus C}=-\varphi_{C},
\end{equation}
whenever $C\in \pud X$.

Straightforward calculations lead to the following
\begin{lemat}
\label{eta}
Let $X$ be a one dimensional box.  Let $A$ and $C$ be elements of $\pud X$. Then
$$
\sum_{B\in \ka O (X)} \eta_B(C)\eta_B (A)= 2^{|X|-2}\iver{A=C}+2^{|X|-3}
\iver{\{A, X\setminus A\} \neq \{C, X\setminus C\}}.
$$
\end{lemat}

This lemma implies that for each $C\in \Pud X$,
\begin{equation}
\label{zwiazek}
\varphi_C= 2^{-|X|+2}\sum_{B\in \ka O(X)}(\eta_B(C)-\eta_B(X\setminus C))\eta_B.
\end{equation}
(Let us emphasize that the functions $\eta_B$  are restricted here to $\pud X$.)
\begin{lemat}
\label{baza1}
Let $X$ be a one dimensional box.
Let $\ka B\subset \Pud X$ be  defined so that for every $A\in \Pud X$, it contains exactly one of the two elements
$A$ and $X\setminus A$. Then the set $H:= \{ \varphi_C: C\in \ka B\}$ is an orthogonal
basis of $\add {X,\er}$.
\end{lemat}

\proof
Since we have already learned that the elements of
$H$ are mutually orthogonal, it remains to show that each $g\in \add {X, \er}$
is a linear combination of them.
As we know from the definition of an additive mapping,
there is a number  $s$ such that $s=g(A)+g(X\setminus A)$  for every $A\in \pud X$.
Let
$$
h= \frac{1}{2}(s \varphi_X +\sum_{C\in \ka B\setminus\{X\}} (g(C)-g(X\setminus C))\varphi_ C).
$$
Fix any $A\in \pud X$. If $A\in \ka B$, then we get
$$
h(A)=\frac{s+ g(A)-g(X\setminus A)}{2}=g(A).
$$
If $A\not\in \ka B$, then by (\ref{sym}), we get
$$
h(A)=\frac{s+ (g(X\setminus A)-g(A))\varphi_{X\setminus A}(A)}{2}=g(A).
$$
Thus,  $g$ coincides with $h$
and consequently $g$ is a linear combination of elements of $H$.
\hfill$\square$

\begin{lemat}
\label{baza2}
Let $X$ be a one dimensional box. Then the set $K:=\{\eta_B\colon \pud X\to \er : B\in \ka O(X)\}$ is a
basis of $\add {X,\er}$.
\end{lemat}

\proof
Let $\ka B$ and $H$ be as defined in Lemma \ref{baza1}. Observe  that
$$
|H|=|K|=2^{|X|-1}.
$$
To complete the proof, it suffices to notice that according to (\ref{zwiazek}), each element of
$H$ is a linear combination of elements of $K$.
\hfill $\square$

\medskip
\medskip
Now, we go back to the general case. Let $\bigotimes_{i=1}^d \add{X_i,\er}$ be
the tensor product of the spaces $\add{X_i, \er}$, where $X_i$, $i\in [d]$,
are the direct factors of $X$. We can  and we do
identify $\bigotimes_{i=1}^d \add {X_i,\er}$ with the subspace of $\er^{\pud X}$ spanned
by the functions $f_1\otimes\cdots\otimes f_d$ defined by
$$
f_1\otimes\cdots\otimes f_d (C)=f_1(C_1)\cdots f_d (C_d),
$$
where $f_i\in \add{X_i, \er}$.

\begin{tw}
\label{struktura}
 Let $X$ be a $d$-box. Then  $\add{X,\er}= \bigotimes_{i=1}^d \add{X_i,\er}$.
\end{tw}
\proof
Let $f\in \add{X,\er}$. Fix $G\in \pud {X_d}$. Let $f_G\colon \pud {X_{d'}}\to \er$ be defined by
$f_G(U)=f(U\times G)$. Clearly, $f_G\in \add{X_{d'}, \er}$. Let $g_1,\ldots , g_m$ be any basis
of $\add{X_{d'}, \er}$. We can write $f_G$ as a linear combination of
the elements of this basis
\begin{equation}
\label{nae}
f_G=\sum_{i=1}^m \alpha_i(G)g_i.
\end{equation}
Let $H\in \pud {X_d}$. Since $f$ is  additive, for each $U$, we have
$$
f_G(U)+f_{X_d\setminus G}(U)=f_H(U)+f_{X_d\setminus H}(U).
$$
Consequently, the mapping $\alpha_i\colon \pud{X_d}\to \er$
is additive for each $i\in [m]$. As
$$
f(U\times G)=\sum_{i=1}^{m} \alpha_i(G)g_i(U)=\sum_{i=1}^{m} g_i\otimes\alpha_i(U\times G),
$$
we obtain  $f\in \add{X_{d'},\er}\otimes\add{X_d,\er}$. Therefore,
$ \add{X,\er}\subseteq\add{X_{d'},\er}\otimes\add{X_d,\er} $. The opposite inclusion is obvious.
Our result follows now by induction with respect to the dimension $d$.
\hfill$\square$

\medskip
\medskip
Let us observe that if $X$ is a $d$-box, then for each $B\in \ka O(X)$, we have
\begin{equation}
\label{eta2}
\eta_B=\eta_{B_1}\otimes\cdots\otimes\eta_{B_d},
\end{equation}
where we interpret $\eta_B$ as defined on $\pud X$. By Lemma \ref{baza2} and the
preceding theorem, we obtain immediately
\begin{tw}
\label{bazad2}
Let $X$ be a $d$-box. Then the set $K:=\{\eta_B: B\in \ka O(X)\}$ is a
basis of $\add {X, \er}$.
\end{tw}

Let $f\in\add{X,\er}$. It follows from our theorem that
there are real numbers $\alpha_B$, $B\in \ka O(X)$, such that
$$
f=\sum_{B\in \ka O(X)} \alpha_B\eta_B.
$$
Since each $\eta_B$ is naturally defined on $2^X$, the above equation determines
the extension $\tilde{f}$ of $f$ to $2^X\setminus\{\emptyset\}$, and the extension $\bar{f}$ of $f$
to the family of all polyboxes.

Let us observe  that by Lemmas \ref{comp} and \ref{symmdiff},
for every polybox $F\subset X$
and every proper suit  $\ka F$ for $F$, we have
$$
\eta_B(F)=\sum_{A\in  \ka F}\eta_B(A),
$$
which implies
\begin{equation}
\label{ekst}
\bar{f}(F)=\sum_{A\in \ka F}f(A).
\end{equation}
In particular, we have the following basic results.
\begin{tw}
\label{add1}
Given a polybox $F$ contained in a $d$-box $X$. For any two proper suits $\ka F$ and $\ka G$
for $F$, and any $f\in \add {X, \er}$,
\begin{equation}
\label{rowno1}
\sum_{A\in \ka F} f(A)=\sum_{A\in \ka G}f(A).
\end{equation}
\end{tw}
\begin{tw}
\label{sdisjoint}
If $F_1,\ldots, F_m$ are pairwise strongly disjoint polyboxes in a $d$-box $X$, then 
$$
\bar{f}(F_1\cup \cdots\cup F_m)=\bar{f}(F_1)+\cdots +\bar{f}(F_m).
$$
\end{tw}  
\begin{tw}
\label{ekst1}
Let $X$ be a $d$-box. Let $\mathbf{A}_P( X,\er)$ be the space of all real-valued functions $g$
defined on the family of all polyboxes contained in $X$ such that for
any polybox $F$, and any proper suit $\ka F$ for $F$
$$
g(F)=\sum_{A\in \ka F}g(A).
$$
Then the  mapping $\add{X, \er}\ni f \mapsto \bar{f}\in \mathbf{A}_P( X,\er)$ is  a linear isomorphism.
\end{tw}

Let us note for completeness that we have a similar theorem for the other extension.

\begin{tw}
\label{ekst2}
Let $X$ be a $d$-box. Let $\mathbf{A}_S(X,\er)$ be the space of all functions
$g: 2^{X}\setminus\{\emptyset\} \to \er$
such that for
any  $F$ and $G\in 2^X\setminus\{\emptyset\}$, if $\widehat{F}$ and $\widehat{G}$ are disjoint,
then
$$
g(F\cup G)=g(F)+g(G).
$$
The mapping $\add{X, \er}\ni f \mapsto \tilde{f}\in \mathbf{A}_S(X,\er)$ is  a linear isomorphism.
\end{tw}

In order to formulate an analogue of Lemma \ref{baza1} for arbitrary dimensions,
we have to introduce some extra terminology.

For a box $C\in \Pud X$ and $\varepsilon \in \zet_2^d$, let
$C^\varepsilon=C^{\varepsilon_1}_1\times\cdots\times C^{\varepsilon_d}_d$ be defined by
$$
C^{\varepsilon_i}_i=\left\{
\begin{array}{cl}
C_i & \text{if $\varepsilon_i=0$,} \\
X_i\setminus C_i & \text{if $\varepsilon_i=1$,}
\end{array}
\right.
$$
where $i\in [d]$. If $C$ is a proper box, then $C^\varepsilon$ is a proper box. If $C$ is not proper, then
it can happen $C^{\varepsilon}$ is empty.
Now, let
$
\varphi_C=\varphi_{C_1}\otimes\cdots\otimes \varphi_{C_d}.
$
Similarly as in the case $d=1$, for any $C$ and $D\in \Pud X$, we have $\varphi_C \perp \varphi_D$,
whenever $C\not\in \{D^\varepsilon : \varepsilon\in \zet_2^d\}$.

For $\varepsilon\in \zet_2^d$, let $|\varepsilon |=|\{i: \varepsilon_i=1\}|$.
The $d$-dimensional counterpart of (\ref{sym}) reads as follows: If
$C\in \Pud X$ and $C^\varepsilon$ is non-empty, then
\begin{equation}
\label{symd}
\varphi_{C^\varepsilon}=(-1)^{|\varepsilon|}\varphi_C.
\end{equation}

\begin{tw}
\label{bazad1}
Let $X$ be a $d$-box.
Let $\ka B\subset \Pud X$ be  defined so that for every $A\in \Pud X$, $\ka B$ contains exactly one element of
the set $\{A^\varepsilon :\varepsilon \in \zet_2^d\}$. Then the set $H:= \{ \varphi_C: C\in \ka B\}$ is an orthogonal
basis of $\add {X,\er}$.
\end{tw}

\proof
For each $i\in [d]$, let us fix  $\ka B_i\subset \Pud {X_i}$ so that for $X=X_i$ and $\ka B=\ka B_i$,
the assumptions of Lemma \ref{baza1} are satisfied.
Then, by Theorem \ref{struktura},
$$
G:=\{\varphi_{C_1}\otimes\cdots \otimes \varphi_{C_d}: C_1\in \ka B_1,\ldots , C_d \in\ka B_d\}
$$
is a basis  of $\add {X, \er}$. On the other hand, one observes that by (\ref{symd})
and the definitions of $H$ and $G$ we have $\varphi_C\in G$ if and only if
exactly one of the two elements $-\varphi_C$, $\varphi_C$ belongs to $H$. Thus, $H$ is a basis.
The orthogonality of $H$ is clear.
\hfill$\square$

\medskip
\medskip
Equation (\ref{zwiazek}) expresses $\varphi_C$ in terms of  functions $\eta_B$ in dimension one.
We apply it to obtain an analogous result which is valid in arbitrary dimensions.
As before, let $X$ be a $d$-box and $C\in\Pud X$.
By (\ref{zwiazek}), (\ref{eta2}) and the definition of $|X|_1$ we have
\begin{eqnarray*}
\varphi_C & = & 2^{-|X|_1+2d}\sum_{B_1\in\ka O(X_1)}\!\!\cdots\!\!\sum_{B_d\in \ka O(X_d)}
\left(\prod_{k=1}^d(\eta_{B_k}(C_k)-\eta_{B_k}(X\setminus C_k))\right)\eta_B\\
          & = & 2^{-|X|_1+2d}\sum_{B\in\ka O(X)}\sum_{\varepsilon\in\zet_2^d}
(-1)^{|\varepsilon|}\eta_B(C^\varepsilon )\eta_B.
\end{eqnarray*}
Hence
\begin{eqnarray}
\varphi_C  =  2^{-|X|_1+2d}\sum_{B\in\ka O(X)}
\sum_{\varepsilon\in\zet_2^d} (-1)^{|\varepsilon|}\eta_B(C^\varepsilon )\eta_B.
\end{eqnarray}

For further use, we introduce
a class of bases of $\add {X,\er}$ containing the bases described in
Theorem \ref{bazad1}.

For each $i$, let $\chi_i\in\add {X_i, \er}$ be the
characteristic function of a non-empty subfamily
$\ka F_i\subseteq \pud {X_i}$. Observe that if $\ka F_i=\pud{X_i}$,
then $\chi_i=1= \varphi_{X_i}$, and that $\chi_i\neq 1$ is equivalent to saying
that, as in Lemma \ref{baza1},
for each $A\in \pud {X_i}$, exactly one of the two elements
$A$, $X_i\setminus A$ belongs to $\ka F_i$. Now, for every
$D\in \Pud {X}$, let
$\tau_D=\tau_{D_1}\otimes\cdots \otimes\tau_{D_d}$, where $\tau_{D_i}=\chi_i$, if $D_i=X_i$,
and $\tau_{D_i}=\varphi_{D_i}$, otherwise.
Clearly, we have
\begin{pr}
\label{bazad3}
Let
$\ka B \subset\Pud {X}$ be defined so that for
every $D\in \Pud{X}$, $\ka B$ has only
one element in common with $\{D^\varepsilon: \varepsilon\in\zet_2^d\}$.
Then the set  $L:=\{\tau_D:  D\in \ka B\}$ is a  basis of $\add {X, \er}$.
\end{pr}

\section{Dyadic labellings}
\label{dilabel}
Let $X$ be a $d$-box and let $E$ be a non-empty set. A mapping $\lambda \colon \pud X \to E$  is said to be a
\textit{dyadic labelling} if it is `onto', and for every $A$, $B$, $C$ and
$D$ belonging to $\pud X$, if $A$, $B$ and $C$, $D$ are twin pairs, and $A\cup
B=C\cup D$, then $\{\lambda(A), \lambda(B)\}=\{
\lambda(C), \lambda(D) \}$.

\begin{pr}
\label{ind3}
Let $X$ be a $d$-box, and $E$ be a set. If
$\lambda \colon  \pud X\to E$
is a dyadic labelling, then for each polybox $F\subseteq X$ and
every  two proper suits $\ka F$, $\ka G$ for $F$, the following equation
is satisfied
$$
\{\lambda(A) \colon  A\in \ka F\}
 =
\{\lambda(A) \colon  A\in \ka G\}.
$$
\end{pr}
\proof
For each $e\in E$, let us define $f_e\colon \pud X\to\er$ as follows
$$
f_e(A)=[e=\lambda(A) ].
$$
It is straightforward from the definition of dyadic labellings 
that $f_e$ is additive. Therefore, by Theorem \ref{add1},
$$
\sum_{A\in \ka F} f_e(A)=
\sum_ {A\in \ka G}f_e(A),
$$
which readily implies our thesis.
\hfill $\square$

\begin{pr}
\label{sekw}
Let $X$ be a $d$-box, and $A$, $B\in\pud X$. If $\lambda:\pud X\to E$ is a dyadic labelling, then there is $D\in \ka C_A$ such that $\lambda (D)=\lambda (B)$ and for each $i\in [d]$, $D_i=B_i$, whenever $B_i\in \{ A_i, X_i\setminus A_i\}$. 
\end{pr}

\proof 
Let $\lambda (B)=e$.
Define a sequence
$B^0=B, B^1,\ldots , B^{d}$ by induction: Suppose that $B^k$
is already defined and $\lambda(B^k)=e$. If $B_{k+1}\in \{A_{k+1}, X_{k+1}\setminus A_{k+1}\}$, then set $B^{k+1}=B^{k}$. If it is not the case,
then 
choose three proper boxes $R$, $S^1$ and $S^2$ so that the following
equations are satisfied
$$
B^k_{(k+1)'}=R_{(k+1)'}=S^j_{(k+1)'},\qquad j=1,2,
$$
$$
R_{k+1}=X_{k+1}\setminus B^{k}_{k+1},\quad
S^1_{k+1}=A_{k+1},\quad
S^2_{k+1}=X_{k+1}
\setminus A_{k+1}.
$$
(Therefore, $B^k$, $R$ and $S^1$, $S^2$ are twin pairs such that
$B^k\cup  R=S^1\cup S^2$.
By the definition of dyadic labellings
$
e\in\{\lambda(S^1), \lambda(S^2) \}.
$) Set
$$
B^{k+1}=\left\{
\begin{array}{ccc}
S^1 & \text{if} & \lambda(S^1)=e, \\
S^2 & \text{if} & \lambda(S^1)\neq e.
\end{array}
\right.
$$
It is clear by the construction that $D=B^d$ has the desired properties.\hfill $\square$  

\medskip
In the sequel, we shall identify each box $A\in \pud X$ with its
coordinates $(A_1,A_2,\ldots,A_d)$. Consequently, $\pud X$ is
identified with
$\prod_{i=1}^{d}\pud {X_i}$,
and  as such can be considered as a $d$-box. Thus, it makes sense
to define $\pudd X:= \pud{\pud X}$. A box $\ka B\in\pudd X$ is
\textit{equicomplementary} if for every $A\in \pud X$ the
intersection $\ka B\cap \ka C_A$ consists of exactly one element.
Equivalently, one can say that for each $i\in [d]$ and each
$C\in\pud {X_i}$ exactly one of the two sets $C$, $X_i\setminus
C$ belongs to ${\ka B}_i$. It is clear that $|\ka B|=(1/2^d)|\pud X|$. Now,
let $\ka F$ be a  partition of the $d$-box $\pud X$ that
consists of  equicomplementary boxes. $\ka F$ has exactly $2^d$ elements,
which means that this partition is minimal. For each $A\in \pud X$,
there is a unique element  $\lambda_{\ka F}(A)\in\ka F$ such that
$A\in \lambda_{\ka F}(A)$. Thus, we have defined the mapping
$\lambda_{\ka F}\colon  \pud X\to {\ka F}$.
\begin{pr}
\label{diada}
The mapping $\lambda_{\ka F}$ is a dyadic labelling.
\end{pr}
\proof
Suppose that $A$, $B\in\pud X$ are dichotomous. Hence
there is $i\in [d]$ for which $A_i=X_i\setminus B_i$. Since
$A_i\in \lambda_{\ka F}(A)_i$, the equicomplementarity of
$\lambda_{\ka F}(A)$ implies that $B_i\not\in \lambda_{\ka
F}(A)_i$. Therefore, $B\not\in \lambda_{\ka F}(A)$, which in turn
implies $\lambda_{\ka F}(A)\neq
\lambda_{\ka F}(B)$. In particular,  the mapping
$\zet_2^d\ni\varepsilon\mapsto\lambda_{\ka F}(A^\varepsilon)$
is one-to-one, and as $|\ka F|=2^d$, our mapping  $\lambda_{\ka F}$
has to be onto.

Suppose now that $A$, $B\in \ka F$ and $C$, $D\in \ka F$ are two twin
pairs such that $A\cup B=C\cup D$. There is an index $i\in [d]$ such 
that 
$$
A_{i'}=B_{i'}=C_{i'}=D_{i'}.
$$
Thus, for $j\in i'$, we have $C_{j}\in\lambda_{\ka F}(A)_{j}$
and $D_{j}\in\lambda_{\ka  F}(A)_{j} $. From the two
elements $C_i$ and, $D_i=X_i\setminus C_i$  one belongs to
$\lambda_{\ka F}(A)_i$. As a result, $C\in\lambda_{\ka F}(A)$
or  $D\in\lambda_{\ka F}(A)$. Consequently,
$$
\lambda_{\ka F}(A)\in\{\lambda_{\ka F}(C),\lambda_{\ka F}(D)\}.
$$
Since  $A$ and $B$,
and $C$ and $D$ play the same role, we may deduce
that
$$
\{\lambda_{\ka F}(A), \lambda_{\ka F}(B)\}=\{\lambda_{\ka F}(C), \lambda_{\ka F}(D)\}.
$$
\hfill $\square$

\begin{tw}
Let $X$ be a $d$-box and let $E$ be a set that consists of $2^d$ elements.
A mapping $\lambda \colon  \pud X\to E$
is a dyadic labelling if and only if  there is a proper suit $\ka F$ for
$\pud X$, composed of equicomplementary boxes, and a bijection $g\colon
\ka F\to E$ such that
$$
\lambda =g\circ\lambda_{\ka F}.
$$
\end{tw}
\proof
The proof of the implication `$\Leftarrow$' is a consequence of Proposition \ref{diada}. In order to show that
the implication `$\Rightarrow$' is  true, it suffices to prove that
$\lambda^{-1}(e)$ is an equicomplementary box for each $e\in E$. 

By Proposition \ref{sekw}, for each $A$, $B\in \pud X$, we have $\lambda (B)\in \lambda(\ka C_A)$. Since $\lambda$ is `onto', $|\lambda (\ka C_A)|=2^d=|\ka C_A|$. Consequently, the set $\lambda^{-1}(e)\cap \ka C_A$ is a singleton.  

Let $G$, $H\in \lambda^{-1}(e)$. Suppose that $L\in\pud X$ and $L_i\in \{G_i, H_i\}$
for each $i\in [d]$. Let $I=\{i: G_i=L_i\}$. By Proposition \ref{sekw} applied to $A=L$ and $B=G$, there is a $K\in \ka C_L$ such that $\lambda(K)=e$ and $K_I=L_I$. Similarly, for $J=\{i: H_i=L_i\}\supseteq [d]\setminus I$, we can show that there is $M\in \ka C_L$ such that $\lambda(M)=e$ and $M_J=L_J$. Since by the previous part $\lambda$ is `one-to-one' on $\ka C_L$ and $I\cup J=[d]$, we have $K=M=L$. Therefore $L\in \lambda^{-1}(e)$, which implies that $\lambda^{-1}(e)$ is a  box. The equicomplementarity of $\lambda^{-1}(e)$ is again a consequence of the previous part of the proof.
\hfill $\square$   

\section{Binary codes}
\label{binary codes}
As before, let $X$ denote a $d$-box. Suppose that for each $i\in [d]$ we have
chosen a mapping
$ \beta_i \colon 2^{X_i}\setminus \{\emptyset , X_i\}\to \zet_2$
such that for each $A_i$,
\begin{equation}
\label{kod}
\beta_i(A_i)+\beta_i({X_i\setminus A_i})=1.
\end{equation}
The mapping $\beta \colon \pud X\to \zet_2^d$ defined by
$$
\beta (A)=(\beta_1 (A_1),\ldots , \beta_d (A_d))
$$
is said to be a \textit{binary code} of  $\pud X$.

Clearly, each binary code is a dyadic labelling. By Proposition \ref{ind3}, if $\ka F$ is  a proper suit
for a polybox $F\subseteq X$, then the set of code words
$$
\beta(\ka F):=\{\beta (A): A\in\ka F\}
$$
depends only on $F$, not on the particular choice of a suit. It yields the following interesting consequence:
\begin{itemize}
\item[]
For $\varepsilon\in \zet_2^d$, let, as before, $|\varepsilon|=|\{i: \varepsilon_i=1\}|$. The number
$|\beta|_k(F):=|\{A\in\ka F: |\beta(A)|=k\}|$ depends only on $F$. In particular, $|\beta|_k(X)={d\choose  k}$.
\end{itemize}
There  are natural examples of binary codes. Assume that a
$d$-box $X$ is of odd cardinality. For each $A\in \pud X$, let us set
$$
\beta_\text{eo}(A)=(|A_1| (\operatorname{mod} 2),\ldots, |A_d| (\operatorname{mod} 2)).
$$
We call  $\beta_\text{eo}$ the \textit{even-odd pattern}. In the case of the even-odd pattern $|\beta|_k(F)$
counts the number of boxes belonging to any suit $\ka F$ for $F$ that have exactly $k$ factors which are the sets of
odd cardinality.

Similarly, we can define the binary code $\beta_\text{ml}$ called the \textit{more-less pattern}
$$
\beta_\text{ml}(A)=([|A_1|>|X_1|/2],\ldots, [|A_d|>|X_d|/2]).
$$
Here $|\beta|_k(F)$ counts the number of boxes belonging to any suit $\ka F$ for $F$ that have exactly $k$ factors
which are the sets having more than a half of elements of the corresponding factors of $X$.

\section{Indices}
\label{indices}

Let a $d$-box $X$, a suit $\ka F\subseteq\pud X$, a
set $F\subseteq X$ and a box $C\subseteq X$ be given. We  define two kinds of
indices: the \textit{index} of $\ka F$
relative to $C$ 
$$
\Index {\ka F}  C=\sum_{A\in \ka F}{\varphi}_C (A),
$$
and the \textit{index} of $F$ relative to $C$
$$ 
\index F C =\tilde{\varphi}_C (F).
$$
Clearly, if $F$ is a polybox and $\ka F$ is a suit for $F$, then $\tilde{\varphi}_C(F)=\bar{\varphi}_C(F)$ and
by (\ref{ekst}) both indices are equal.

Let us set
$$
\ka C= \ka C_C:=\{C^\varepsilon: \varepsilon\in\zet_2^d\}\setminus\{\emptyset\}.
$$
Obviously, $\ka C$ is a suit for $X$. We  call it \textit{simple}.
Observe that $\ka C$ is proper if and only if $C$ is  a proper box. We say that elements
$A$ and $B$ of $\ka C$ \textit{carry the same sign} if
$(-1)^{|\varepsilon|}=1$ for the only $\varepsilon$ such that $A^\varepsilon=B$.
If $C\neq X$, then \textit{carrying the same sign} is an equivalence relation with two classes
of abstraction. Any assignment to one of these classes
$+1$ and to the other $-1$ is called an \textit{orientation} of $\ka C$.
Let $\ka C^+$ be the class to which $+1$ is assigned
and $\ka C^-$ be the class to which $-1$ is assigned.
Suppose now that the orientation of $\ka C$ is \textit{induced} by $C$, that is, $C\in \ka C^+$.
Let $$
I=i(C):=\{i\in[d]: C_i\neq X_i\}.$$
Let
$ n^+(\ka F, C)$ be the number of all those $A\in \ka F$ that
$A_I\in  (\ka C^+)_I$
and accordingly,  $n^-(\ka F, C)$ be the number of all those $A\in \ka F$ that $A_I\in  (\ka C^-)_I$.
Immediately from the definition of the indices we obtain
\begin{pr}
\label{ladunek1}
Let $X$ be a $d$-box and let $\ka F$ be a proper suit for a polybox $F$. For any $C\in \Pud X$,
if $C\neq X$, then
\begin{equation}
\label{ladunek2}
\index F C =\Index {\ka F} C= n^+(\ka F, C)-n^-(\ka F, C).
\end{equation}
Moreover,
$$
\index F X=|\ka F|=|F|_0.
$$
\end{pr}

In the case of $C\in \pud X$ the meaning of (\ref{ladunek2}) is particularly
simple, it counts the sum of all signs that come from those elements of $\ka C$
that belong to the  suit $\ka F$ for $F$.
In particular, this sum is independent of the choice of a suit for $F$.

\begin{wn}
\label{zerot}
For each proper suit $\ka F$ for a $d$-box $X$ and each $C\in \Pud X$, if $C\neq X$, then
$\Index {\ka F} C=0$. In particular, if $C$ is a proper box
then  the size of ${\ka F}\cap {\ka C_C}$ is even.
\end{wn}






\begin{lemat}
\label{2kostki}
If $F$ and $G$ are two polyboxes in a $d$-box  $X$, and $\{\zeta_i:
i\in [2^{|X|-d}]\}$ is a basis of $\add {X, \er}$, then
$F=G$ if and only if $\bar{\zeta}_i(F)=\bar{\zeta}_i(G)$ for each $i\in[2^{|X|-d}]$.
\end{lemat}

\proof Let $x:=(x_1,\ldots ,x_d)$ be an arbitrary element of $X$, and let
$\bar{x}=\{x_1\}\times\cdots\times\{x_d\}$. Then
$\bar{x} \in \ka O(X)$ and $\eta_{\bar{x}}\in \add {X,\er}$.
Obviously, we have $[x\in F]=\eta_{\bar{x}}(F).$
Since $\zeta_i$, $i\in [2^{|X|-d}]$, is a basis, there are reals $\alpha_i(x)$ such that
$\eta_{\bar{x}}=\sum_i \alpha_i(x)\zeta_i$. Thus,
$$
[x\in F]=\sum_i \alpha_i(x)\bar{\zeta}_i(F).
$$
Our result follows now immediately, as the same equation holds true for $G$. \hfill $\square$
\begin{tw}
\label{kryt1}
If $F$ and $G$ are two polyboxes in a $d$-box  $X$, then
$F=G$ if and only if $\index FC=\index GC$ for each $C\in \Pud X$.
\end{tw}

\proof
Let $\ka B\subset \Pud X$ be as in Theorem \ref{bazad1}. Recall the functions $\varphi_C$, $C\in \ka B$, form
an (orthogonal) basis of $\add{X,\er}$, and for each $C\in \Pud X$,
$$
\index F C =\bar{\varphi}_C (F).
$$
Therefore, our theorem is a consequence of the preceding lemma.\hfill $\square$

\section{Symmetric polyboxes}
\label{symbox}
Let $X$ and $Y$ be $d$-boxes. Suppose that for each $i\in [d]$  one has given a mapping $\varphi_i\colon\pud {X_i}\to\pud{Y_i}$ such that
$$
\varphi_i(X_i\setminus A)=Y_i\setminus \varphi_i(A),
$$
whenever $A\in \pud {X_i}$. Then one can define the mapping $\Phi:=\varphi_1\times\cdots\times\varphi_d$ from $\pud X$ into $\pud Y$:  
$$
\Phi (B)=\varphi_1(B_1)\times\cdots\times \varphi_d(B_d). 
$$
It is clear that $\Phi$ sends  dichotomous boxes into  dichotomous boxes. Moreover,  if $A$, $B$ and $C$, $D$  are two twin pairs in $\pud X$ such that $A\cup B = C\cup D$, then $\Phi(A)$, $\Phi(B)$ and $\Phi(C)$, $\Phi(D)$ are two twin pairs as well, and $\Phi(A)\cup\Phi(B)=\Phi(C)\cup\Phi(D)$.  Further we shall refer to $\Phi$ as a mapping that \textit{preserves dichotomies}. 

\begin{pr} 
\label{rowy}
Given two $d$-boxes $X$, $Y$ and a mapping $\Phi\colon \pud X\to \pud Y$, that preserves dichotomies.   
If  $\ka F$, $\ka G\subset \pud X$ are proper suits and $\bigcup \ka F=\bigcup \ka G$, then their images 
are also proper suits and $\bigcup \Phi(\ka F)=\bigcup\Phi(\ka G)$. 
\end{pr}
\proof
The proof is much the same as that of Proposition \ref{ind3}.  For each $y\in Y$, let us define $e_y\colon\pud X\to \er$ by the formula $e_y(A)=[y\in \Phi (A)]$. Since $\Phi$ preserves dichotomies,  $e_y$ is additive. By Theorem \ref{add1},
$$
\sum_{A\in \ka F} e_y(A)=\sum_{A\in \ka G}e_y(A).
$$ 
Hence $\bigcup\Phi(\ka F)=\bigcup \Phi(\ka G)$.\hfill $\square$ 

\medskip
By the above proposition, we can extend $\Phi$ to polyboxes. Namely, if $F$ is a polybox in $X$ and $\ka F$ is a suit for $F$, then we may define $\Phi(F)$ by the equation
$$
\Phi(F)=\bigcup \Phi(\ka F).
$$
In particular, if $\varepsilon\in \zet_2^d$, then the polybox $F\naw \varepsilon\subseteq X$ given by
$$
F^{(\varepsilon)}=\bigcup_{A\in \ka F}A^\varepsilon
$$
is well-defined. Observe that if $F$ is an improper box, then $F^{(\varepsilon)}$ does not have to coincide with already defined $F^\varepsilon$, as the latter set is empty  if there is an $i\in [d]$ such that $F_i=X_i$ and $\varepsilon_i=1$, while $F\naw \varepsilon$ is non-empty.

\begin{pr}
\label{commute0}
Given two $d$-boxes $X$, $Y$ and a mapping $\Phi\colon \pud X\to \pud Y$ that preserves dichotomies. Then for every polybox $F\subseteq X$ and $\varepsilon\in \zet_2^d$
\begin{equation}
\label{commute}
\Phi(F^{(\varepsilon)})=\Phi(F)^{(\varepsilon)}.
\end{equation}
\end{pr}
\proof It suffices to observe that (\ref{commute}) is satisfied when $F$ is a proper box.\hfill $\square$   

\begin{lemat}
\label{suma00}
Let $U$ be an $n$-box, $f\colon U\to \zet_2$ and $V=\{u\in U: f(u)=1\}$. If for every $x\in U$ and every $i\in[n]$
\begin{equation}
\label{suma0}
\sum_{u\in U} f(u)[u_{i'}=x_{i'}]=0,
\end{equation}    
and $V\neq \emptyset$, then $|V|\ge 2^n$. Moreover, for every $v\in V$ there is some $y\in V$ such that 
$v_i\neq y_i$, $i=1,\ldots,  n$.   
\end{lemat}
\proof
Let us fix $v\in V$, and define $\gamma\colon V\to \{0,1\}^n$ by the formula
$$
\gamma(u)_i=[u_i\neq v_i].
$$ 
It suffices to show that $\gamma$ is `onto'. Suppose that it is not. Then there is an element $\tau$ with minimal support among the members of $\{0,1\}^n\setminus \gamma (V)$. Let us choose an index $i$ such that $\tau_i=1$ (such $i$ has to exist as $\gamma(v)=0$). Let $\tau'$ be defined so that $\tau'_i=0$ and $\tau'_{i'}=\tau_{i'}$. Since $\tau'$ has its support strictly contained in that of $\tau$, it belongs to $\gamma(V)$. Therefore, there is an $x\in V$ such that $\tau'=\gamma(x)$.
By (\ref{suma0}), there has to exist $u\in V\setminus\{x\}$ for which $u_{i'}=x_{i'}$. By the definition of $\gamma$, and that of $\tau'$, we conclude that $\gamma(u)=\tau$, which is a contradiction. \hfill $\square$      

\begin{tw}
\label{wieza}
Let $F$ be a polybox in a $d$-box $X$ and let $\ka B\subset \Pud X$ intersect with each $\ka C_C$, $C\in \Pud X$, at exactly one element. Let $I\subseteq [d]$ be a set of size $n\ge 2$, and
$$
\ka E=\{C\in \ka B: i(C)=I,\, \index F C \equiv 1 \mod 2\}.
$$ 
Suppose that for each $B\in \ka B$ 
$$
\index F B\equiv 0 \mod 2,
$$
whenever $i(B)\subset I$ and $|i(B)|=n-1$. Then
$\ka E$ is empty or $|\ka E|\ge 2^n$. Moreover, if $\ka E$ is non-empty, then $n<d$ and for every $C\in \ka E$ there is a $D\in \ka E$ such that  $D_i\not\in\{C_i, X_i\setminus C_i\}$ for each $i\in I$.
\end{tw}

\proof 
For each $i\in[d]$ let us pick ${\ka B}_i\subseteq \Pud {X_i}$ so that $\{A, X_i\setminus A\}\cap {\ka B}_i$ is a singleton for each $A\in \Pud {X_i}$. It is clear that it suffices to  proceed with  $\ka B$ defined as follows
$$
\ka B=\{ B\in \Pud X: B_i\in \ka B_i\,\,  \text{for each $i\in [d]$}\},
$$
as for every $B$ and $C\in \Pud X$ if $B\in \ka C_C$, then 
$$
\varphi_B\equiv \varphi_C \mod 2.
$$
(Observe that  $\ka B$ can be thought of as an equicomplementary box, which was defined in Section \ref{dilabel}.)   
Let us order the elements of $I$: $i_1 < \ldots < i_n$. Let $\ka B^\circ=\{C\in \ka B: i(C)=I\}$ and     
$
U=\{(C_{i_1},\ldots , C_{i_n}): C\in \ka B^\circ \}.
$
Define $f\colon U\to \zet_2$ by the formula 
$$
f(C_{i_1},\ldots , C_{i_n})=\index F C \przyst 2,
$$
where $C\in \ka B^\circ$. Let $V=\{(C_{i_1},\ldots ,C_{i_n}): C\in \ka E\}$.  

It is clear that $f$ just defined satisfies the assumptions of Lemma \ref{suma00} if and only if the following claim holds true
\vskip 1em

\noindent\textbf{Claim}\,\,
For every $C\in \ka B^\circ$ and $i\in I$
\begin{equation}
L:=\sum_{D\in \ka B^\circ}\index F D [D_{i'}=C_{i'}]\equiv 0 \mod 2.
\end{equation}  
\vskip 1em

To prove it, let us define $B\in \ka B$ so that $B_i=X_i$ and $B_{i'}=C_{i'}$.
For $D\in \Pud X$, let $\ka D_D$ be defined by the following equivalence 
$$
A\in \ka D_D \Leftrightarrow A\in \pud X\,\,\text{and there is a $E\in \ka C_D$ such that  $E_{i(D)}=A_{i(D)}$.}
$$ 

Let $\ka F$ be a proper suit for $F$.
It follows from Proposition \ref{ladunek1} that
$$
\index F D\equiv |\ka F\cap \ka D_D |\mod 2.
$$  
Moreover,  since the sets $\ka D_D$, $D\in \ka B^\circ$, are pairwise disjoint and 
$$
\bigcup\{\ka D_D: D\in \ka B^\circ,  D_{i'}=C_{i'}\}=\ka D_B,
$$
we get
$$
L\equiv\sum_{D\in \ka B^\circ} |\ka F\cap \ka D_D |[D_{i'}=C_{i'}]=|\ka F\cap\ka D_B|\equiv \index F B \mod 2.
$$ 
But the latter number is even by the assumption.

Now, from Lemma \ref{suma00} we conclude that $\ka E$ is either empty or has at least $2^n$ elements. If it could happen that $\ka E$ is non-empty while $n=d$, then, again by Lemma \ref{suma00}, for any $G\in\ka E$ there would exist another box $H\in \ka E$ such that $G_i\neq H_i$ for each $i\in [d]$. This fact and the definition of $\ka E$ would imply that 
$\ka C_G$  and $\ka C_H$ would intersect $\ka F$, which would be  impossible as if $A\in \ka C_G$ and $B\in \ka C_H$, then $A$ and $B$ are not dichotomous.\hfill $\square$ 

\begin{tw}
Let $F$ be a polybox in a $d$-box $X$ and let  $\varepsilon\in\zet_2^d$ be different from $0$. If $F^{(\varepsilon)}=F$,
then  $\index F B\equiv 0 \mod 2$ for every box $B\in \Pud X\setminus \{X\}$. If in addition the number $n:=i(B)\cap \supp \varepsilon$ is odd, then $\index F B=0$.  
\end{tw}   

\proof
First we prove the second part of the theorem. 

Let $\ka F$ be a proper suit for $F$. Since $F^{(\varepsilon)}=F$, we obtain
$$
\index F B =\sum_{A\in \ka F}\varphi_B(A)=\sum_{A\in \ka F}\varphi_B(A^\varepsilon).
$$
It  follows from the definitions of $\varphi_B$ and $n$ that 
$$
\varphi_B(A^\varepsilon)=(-1)^n\varphi_B(A)=-\varphi_B(A). 
$$ 
Therefore, $\index F B=-\index F B$ and consequently $\index F B=0$.

To prove the first part, we begin with showing that $\index F B\equiv 0 \mod 2$ for any box $B$ for which $i(B)$
is a proper subset of $[d]$. To this end, let us fix $i\not\in i(B)$, and define a $d$-box $Y$ by the equations $Y_{i'}=X_{i'}$, 
$Y_i=\{0,1\}$. Now, choose $\ka B_i$ as in the proof of the preceding theorem and define 
$\Phi\colon \pud X\to \pud Y$ so that if $B=\Phi(A)$, then $B_{i'}=A_{i'}$ and 
$$
B_i=\left\{\begin{array}{ll}
\{0\} & \text{if $A_i\in \ka B_i$,}\\
\{1\} & \text{otherwise.}      
\end{array}\right.
$$
Clearly, $\Phi$ preserves dichotomies. Let $C=\Phi(B)$. As is seen from the definition of $\Phi$, for each 
$A\in \pud X$, one has $\varphi_B(A)=\varphi_C(\Phi(A))$. Therefore, if $G=\Phi (F)$, then
\begin{equation}
\label{zapis}
\index F B=\index G C.
\end{equation}  
Let $G_0=\{y\in G: y_i=0\}$ and $G_1=\{y\in G: y_i=1\}$. These two sets are polyboxes.  
By Proposition \ref{commute0} and the fact that $F^{(\varepsilon )}=F$, we obtain $G^{(\varepsilon)}=G$. If 
$\varepsilon_i=0$, then 
\begin{equation}
\label{warstwy}
G_{\kappa}^{(\varepsilon )}=G_{\kappa},\,\, \text{for $\kappa=0,1$.}
\end{equation}
If we ignore the $i$-th coordinate, which is constant for elements of both polyboxes $G_\kappa$, then we can think of them as polyboxes in a $(d-1)$-box. Thus, by induction and (\ref{warstwy}), we can maintain that 
$\index {G_\kappa} C \equiv 0 \mod 2$. Consequently, this,  together with (\ref{zapis}), yields           
\begin{equation}
\label{dublet}
\index F B =\index G C =\index {G_0} C +\index {G_1} C\equiv 0 \mod 2. 
\end{equation}
If $\varepsilon_i=1$, then $G^{(\varepsilon)}=G$ implies $G_0^{(\varepsilon)}=G_1$ and $G_1^{(\varepsilon)}=G_0$.
Hence
$$
\index {G_0} C\equiv \index {G_1} C \mod 2 
$$
and (\ref{dublet}) holds true also in this case. 

The fact that $\index F B\equiv 0 \mod 2$ when $i(B)=[d]$ follows from what has been proved before and Theorem~\ref{wieza}.\hfill$\square$     
\section{Rigidity}
\label{rikitiki}
\begin{lemat}
\label{znaki}
Let $\ka F\subset\Pud X\setminus\{X\}$ be a suit for a $d$-box $X$ and let $C\in \ka F$ be chosen so that 
$$
k:=|i(C)|=\max\{|i(D)|\colon D\in \ka F\}.
$$
For each $D\in \ka C_C$, let $\varepsilon (D)\in \zet_2^d$ be defined by the equation $D=C^{\varepsilon(D)}$.
Then 
$$
\sum_{D\in \ka C_C\cap\ka F} (-1)^{|\varepsilon(D)|}=0.
$$
\end{lemat}
\proof
Observe first that for any $D\in \Pud X$,  $\index D C=0$ if and only if $i(D)\not\supseteq i(C)$. If $D\in \ka F$ and $i(D)\supseteq i(C)$, then by the definition of $k$, $D\in \ka C_C$. In this case, by Proposition \ref{ladunek1},  $\index D C =(-1)^{|\varepsilon(D)|}2^{d-k}$.
Since $C\neq X$, it follows from Proposition \ref{zerot} that $\index X C=0$.  Consequently, 
by Theorem \ref{sdisjoint}, we obtain
$$
0=\sum_{D\in \ka F} \index D C=
2^{d-k}\sum_{D\in \ka C_C\cap \ka F}  (-1)^{|{\varepsilon (D)}|}.
$$ 
\hfill $\square$

\begin{tw}[rigidity]
\label{sztywnypal}
Let $F$ be a polybox in a $d$-box $X$, and $\ka F$ be a suit for $F$. Suppose that 
$$
|\index F C|=|\ka F\cap \ka C_C|,
$$
for each $C\in \pud X$. Then $\ka F$ is the only suit for $F$.
\end{tw}

Let $X$ be a one dimensional box. Let $\ka E(X)\subseteq \pudd X $ be the family of all equicomplementary boxes. For every $A\in \pud X$, we define  $\ka E A$ as the set of all these equicomplementary boxes $B\in \ka E(X)$ that $A\in B$. 

\begin{lemat}
\label{intermediolan}
If $X$ is a one dimensional box, $A^1,\dots,A^n$ belong to $\pud X$ and  $A^i\not\in\{A^j, X\setminus A^j\}$ for every $i\neq j$, then  
$$
|\ka E A^1 \cap\dots\cap \ka E A^n|=\frac{1}{2^{n}}|\ka E(X)|.
$$    
\end{lemat}
\proof
Let $\ka A=\{ \{A, X\setminus A\}\colon A\in \pud X\}$. In order to form an equicomplementary box  we pick independently one  element from each set belonging to $\ka A$. Let $x_1,\dots x_n$ be different elements  of $\ka B$. If for each $i\in [n]$ we fix one element of $x_i$, then the number of equicomplementary boxes which can be formed so that they contain these fixed elements is $\frac{1}{2^{n}}|\ka E(X)|$. Clearly, by our assumptions, the elements $x_i=\{A^i, X\setminus A^i\}$, $i\in [n]$, are different. \hfill $\square$          

\medskip
If $X$ is  a $d$-box and $A\in\Pud X$, then we define $\breve{A}=\ka E A_1\times\cdots\times \ka E A_d$.

\begin{lemat}
\label{baran}
Let $X$ be a $d$-box and let $B$, $C$, $D\in \pud X$ be such that $C\neq D$, and $B$, $D$ are not dichotomous. Then 
$$
\breve{B}\cap\breve{C}\neq \breve{B}\cap\breve{D}.
$$
\end{lemat}     

\proof If it were true that $\breve{B}\cap\breve{C}= \breve{B}\cap\breve{D}$, then 
$$
|\breve{B}\cap\breve{C}\cap\breve D|= |\breve{B}\cap\breve{D}|.
$$
Therefore,
$$
\prod^d_{i=1}|\ka E B_i\cap\ka E C_i\cap \ka E D_i|=\prod^d_{i=1}|\ka E B_i\cap \ka  E D_i|,
$$
As $B$ and $D$ are not dichotomous, the latter product is different from zero. Thus, for each $i$,
$$ 
|\ka E B_i\cap\ka E C_i\cap \ka E D_i|=|\ka E B_i\cap \ka  E D_i|,
$$
Now, it follows from the preceding lemma that $\ka E C_i=\ka E D_i$, for each $i$, which readily implies $C=D$.
\hfill$\square$

\medskip
\dowod {of Theorem \ref{sztywnypal}} The mapping $\Phi\colon \pud X\to \pud {\breve{X}}$  given by 
$
\Phi (A)=\breve{A}
$
preserves dichotomies. Suppose that there are two different suits $\ka F$ and $\ka G$ for $F$. Then there is a $B\in \ka G\setminus \ka F$. Moreover, since by Proposition \ref{rowy} $\bigcup \Phi(\ka G)=\bigcup \Phi (\ka F)$, we deduce that $\breve{B}\subset \bigcup \Phi(\ka F)$. The set $\breve{B}$ is a $d$-box and $\ka H=\{\breve{B}\cap\breve {A}\colon A\in \ka F\}\setminus \{\emptyset\}$ is a (possibly improper) suit for $\breve{B}$. Clearly, 
$\breve{B}\not\in \ka H$. By Lemma \ref{znaki}, there is a $C\in \ka H$ such that 
\begin{equation}
\label{dupa}
\sum_{D\in \ka C_C\cap \ka H} (-1)^{|\varepsilon (D)|}=0,
\end{equation}       
where we interpret $\ka C_C$ as a subfamily of $\Pud {\breve{B}}$. In accordance with this interpretation,
for every $\varepsilon\in \zet_2^d$ and each $i\in [d]$,
$$
(C^\varepsilon)_i=
\left\{
\begin{array}{ll}
C_i & \textrm{if $\varepsilon_i=0$,}\\
\breve{B}_i\setminus C_i & \textrm{if $\varepsilon_i=1$.}
\end{array}
\right.
$$   
Let $U\in \ka F$ be chosen so that $C=\breve{B}\cap \breve{U}$. Observe that 
$$
D=C^{\varepsilon (D)}= \breve{B}\cap (\breve{U})^{\varepsilon(D)},
$$
for each $D\in \ka C_C\cap \ka H$.
As $\Phi$ preserves dichotomies, by Proposition \ref{commute0}, we get
$D=\breve{B}\cap \Phi(U^{\varepsilon(D)})$. Since $D\in \ka H$, there is a $K\in \ka F$ for which one has
$D=\breve{B}\cap\breve{K}$. From Lemma \ref{baran} it follows that $K=U^{\varepsilon (D)}$. Hence $\{U^{\varepsilon (D)}:D\in \ka C_C\cap \ka H\}\subseteq \ka F$. Consequently, by (\ref{dupa}) and the definition of the index, we deduce that 
$$
|\index F U|=|\Index {\ka F} U |<|\ka C_U\cap \ka F|.
$$    
\hfill $\square$

\section{Modules}
\label{modules}
Given a non-empty set $Y$. We denote by  $\zet Y$ the free $\zet$-module
generated by $Y$. Every function $f\in \er^Y$ extends uniquely to the
$\zet$-linear function $f^\dag\colon\zet Y\to \er$. Let now $X$ be a $d$-box.
Let $E,F\in \pud X$ and $G,H\in \pud X$ be any two twin pairs such that $E\cup F= G\cup H$.
If $f\in \add{X,\er}$, then by the definition of additive mappings one has
$
f^\dag(E+F-G-H)=0.
$
Let $n(X)$ be the submodule of $\modgen X$ generated by the elements $E+F-G-H$,
where $E,F$ and $G,H$ are as described above.
Clearly, for each $f\in \add{X,\er}$,
$n(X)$ is contained in the set of zeros of $f^\dag$. We show that $n(X)$
is the set of common zeros of the functions $f^\dag$, $f\in \add{X,\er}$:
\begin{pr}
\label{zera}
For every $x\in \modgen X\setminus n(X)$  there is an $f\in \add{X,\er}$ such that $f^\dag(x)\neq 0$.
\end{pr}

\proof
Suppose it is not true. Then one could find the smallest integer $d$ for which there are a $d$-box $X$,
an element $x\in \modgen X\setminus n(X)$ which is a common zero of the functions $f^\dag$, $f\in \add{X,\er}$.
Fix a one dimensional equicomplementary box $\ka B\subset \pud {X_d}$ and $D\in \ka B$.
The element $x$ can be expressed as follows
$$
x=\sum_{E\colon E_d\in \ka B} k_EE+\sum_{F\colon F_d\not\in \ka B} k_FF,
$$
where $E, F\in \pud X$ and $k_E, k_F$ are integers. Observe that for each $F$ one has
$$
F\equiv F_{d'}\times D + F_{d'}\times (X_d\setminus D)- F_{d'}\times(X_d\setminus F_d)\mod n,
$$
where $n=n(X)$.
If we now replace each $F$ by the expression on the right, then the formula for $x$ can be rewritten as follows
$$
x\equiv\sum_{E\colon E_d\in \ka B} l_EE+\sum_{G\in \pud {X_{d'}}} m_G(G\times D +G\times (X_d\setminus D)) \mod n,
$$
where $l_E, m_G$ are properly chosen integers. Fix an element
$A\in \ka B\setminus\{D\}$, then for each $g\in \add{X_{d'},\er}$  we get
$$
0=(g\otimes \varphi_{A})^\dag(x)=\sum_{E\colon E_d=A} l_E g(E_{d'}),
$$
where $\varphi_A$ is as defined in (\ref{fione}).
Let $u=\sum_{E\colon E_d=A} l_E E_{d'}$; therefore, $u$ is a common zero of all mappings $g^\dag$,
$g\in \add{X_{d'},\er}$. Now, it follows from the assumption on the minimality of $d$ that $u\in n(X_{d'})$. Hence
$$
u\otimes A=\sum_{E\colon E_d=A}l_EE\in n,
$$
and to each $G\in\pud {X_{d'}}$ there correspond integers $p_G$, $r_G$ such that
$$
x\equiv  \sum_{G\in \pud {X_{d'}}} p_G(G\times D) + r_G(G\times (X_d\setminus D))\mod n.
$$
Denote the element on the right by $y$ and again take $g\in \add{X_{d'}, \er}$. We have
$$
0= (g\otimes \chi_{\ka B})^\dag(y)=\sum_{G\in \pud {X_{d'}}} p_Gg(G).
$$
By the same argument as before, we deduce that $\sum_{G\in {\pud {X_{d'}}}} p_G(G\times D)\in n$.
Similarly,
$$
0= (g\otimes(1- \chi_{\ka B}))^\dag(y)=\sum_{G\in {\pud {X_{d'}}}} r_Gg(G),
$$
and $\sum_{G\in \pud {X_{d'}}} r_G(G\times(X_d\setminus D))\in n$.
Hence $x\in n$, which is a contradiction. 

\hfill\ $\square$

Let $\addb {X, \er}$ consists of all functions $f\colon\Pud X\to \er$ such that for every twin pair $A$, $B\in \Pud X$,
$$
f(A\cup B)=f(A)+f(B).
$$
\begin{lemat}
\label{equiv}
If $f\in \addb {X,\er}$, then there is a function $g\in\add {X,\er }$ such that $f=\bar{g}|\Pud X$.
\end{lemat}
\proof
Let $C\in \Pud X$. If $C$ is not proper, then there is an index $i\not\in i(C)$. Let us choose a twin pair
$A$, $B\in \Pud X$ so that $C_{i'}=B_{i'}=A_{i'}$. Then by the fact that $f\in \addb {X,\er }$ we have
$f(C)=f(A)+f(B).$
Observe that $i(A)$ and $i(B)$ are both of greater size than $i(C)$.  If  $A$ and  $B$ are not already
proper boxes, then we split up each of them into a twin pair. We proceed this way  till  we obtain a
proper suit $\ka F$ of $C$.
Clearly, for this suit we have
$$
f(C)=\sum_{F\in \ka F} f(F),
$$
Since $g=f|\pud X$ belongs to $\add {X,\er}$, the latter expression for $f(C)$ shows that
$f=\bar{g}|\Pud X$.\hfill $\square$

In the remainder of this paper, any element from $\add{X,\er}$ and the corresponding element from $\addb{X,\er}$
will be denoted by the same symbol.

Let us denote by $N(X)$ the submodule of $\Modgen X$ spanned by the elements
of the form $(A\cup B)-A-B$, where $A,B\in \Pud X$ run over all twin pairs.

\begin{lemat}
\label{modulo}
$\quad\modgen X+ N(X)=\Modgen X.$ 
\end{lemat}
\proof
It suffices to prove that for each $C\in \Pud X$ there is a $y\in \modgen X$ such that $C-y\in N(X)$. 
We proceed  by induction with respect to $k=d-i(C)$. Let us choose a twin pair $A$, $B\in \Pud X$ so that $C=A\cup B$. As $i(A)$ and $i(B)$ are equal to $i(C)+1$, it follows by the induction hypothesis that there are $y'$ and $y''$ in $\modgen X$ such that $A-y'$ and $B-y''$ are both in $N(X)$. Since in addition $C-A-B\in N(X)$, we obtain that
$$ 
C-(y'+y'')\in N(X),
$$
which completes the proof.\hfill $\square$

\begin{tw}
\label{rzut}
Let $\ka F\subset \pud X$ be an equicomplementary box. Let $\ka B\subset \Pud X$ be
defined so that $\ka B_i=\ka F_i\cup\{X_i\}$, whenever $i\in[d]$. For any $D\in \ka B$,
let $\tau_D\in \add {X, \er}$ be defined so that
$\tau_{D}=\tau_{D_1}\otimes\cdots\otimes\tau_{D_d}$ and
$$
\tau_{D_i}=[D_{i}\in \ka F_i]\varphi_{D_i}+[D_i=X_i](1-\chi_{\ka F_i}),
$$
whenever $i\in[d]$.
Then $\Modgen X=N(X)\oplus \zet\ka B$.
Moreover, the mapping
$P\colon \Modgen X\to \zet\ka B$ given by the formula
$$
P(x)=\sum_{D\in \ka B}{\tau}^\dag_D (x)D
$$
is a projection onto $\zet\ka B$, $\ker P=N(X)$ and  $P(\modgen X)=\zet\ka B$.
\end{tw}
\proof
Since ${\tau}^\dag_D$, $D\in \ka B$, are $\zet$-linear, $P$ is $\zet$-linear as well.
Thus, it suffices to show that: $(1)$ $P(A)=A$ for $A\in \ka B$; $(2)$ $\ker P=N(X)$; $(3)$ $P(\modgen X)=\zet\ka B$.
In order to prove $(1)$, observe first that this claim  is easily seen for $d=1$. Therefore, we may assume that $d>1$.
By the definition of $P$, one can write
$$
P(A)=\sum_D{\tau}_{D_{d'}}(A_{d'}){\tau}_{D_{d}}(A_d)A_{d'}\times A_d=u\otimes v,
$$
where $
u=\sum_{E\in\ka B_{d'}}{\tau}_E(A_{d'})A_{d'}$,  $v= \sum_{F\in \ka B_{d}}{\tau}_F(A_{d})A_d$. By induction,
$u\in \zet\ka B_{d'}$, while $v\in \zet\ka B_d$. Consequently, $P(A)\in \zet\ka B$.

From Lemma \ref{modulo} it follows that for each $x\in \Modgen X$ there is a $y\in \modgen X$ such
that $x\equiv y (\operatorname{mod} N(X))$. Thus, by the definition of $P$, we have
$P(x)=P(y)$, which proves $(3)$.

To prove $(2)$ assume in addition that $x\in \ker P$. Then $y\in \ker P$ as well. By the definition of $P$,
it means that $\tau_D^\dag(y)=0$ for each $D\in \ka B$. By Lemma \ref{bazad3}, the system $\tau_D$, $D\in \ka B$,
is a basis of $\add {X, \er}$. Thus, combining this fact with Proposition \ref{zera} leads to the conclusion
$y\in n(X)$. Since $n(X)\subset N(X)$, we obtain $x\in N(X)$.
\hfill $\square$

\begin{uwa}
\label{kryt2}
Let $\ka F$, $\ka G\subset \pud X$ be two suits. Theorem \ref{kryt1} can be rephrased
so that it can be used to decide whether $\ka F$ and $\ka G$ are suits for the same polybox:
it is the case if and only if $\Index {\ka F} C=\Index {\ka G} C$, for  every $C\in \Pud X$.
This criterion relates to the system $\varphi_C$,  $C\in \Pud X$. As indicated by Lemma \ref{2kostki},
we can replace the latter system by any  basis of $\add {X, \er}$, for example, by $\tau_D$, $D\in \ka B$,
where $\ka B$ is as in Theorem \ref{rzut}. In fact, this theorem gives us a simple computational method to verify
whether two suits determine the same polybox. To describe this method, we shall identify, as we tacitly have
already done,
$\bigotimes_{i=1}^d \Modgen{{X_i}}$ with $\Modgen X$, by letting  $B\in\Pud X$ correspond to
$B_1\otimes\cdots\otimes B_d$, and extending this correspondence by linearity.
Let $P_i: \Modgen{{X_i}}\to \zet\ka B_i$ be the projection with the kernel $N(X_i)$, for each $i\in [d]$.
It can be seen that $P=P_1\otimes\cdots\otimes P_d$ (see the next remark). By Theorem \ref{rzut} and Lemma \ref{2kostki},
two, possibly improper, suits $\ka F$ and $\ka G$ define the same polybox if and only if
$$
\sum_{A\in \ka F} P_1(A_1)\otimes\cdots\otimes P_d(A_d)=\sum_{A\in \ka G} P_1(A_1)\otimes\cdots\otimes P_d(A_d).
$$
Observe that $P_i(A_i)$, $i\in[d]$, are easily calculated
$$
P_i(A_i)=
\begin{cases}
A_i & \text{if $A_i\in {\ka B}_i$},\\
X_i-(X_i\setminus A_i) & \text{if $A_i\not\in {\ka B_i}$}.
\end{cases}
$$
Therefore, in order to decide if $\ka F$  and $\ka G$ are suits for the same polybox it suffices
to write the sum $\sum_{A\in \ka F} P_1(A_1)\otimes\cdots\otimes P_d(A_d)$, where $P_i(A_i)$ are as described above,
and expand it to get an expression of the form $\sum_{C\in \ka B}\alpha_CC$, then repeat the same for $\ka G$.
If the resulting expressions coincide, then the suits define the same polybox, otherwise they do not.
\end{uwa}
\begin{uwa}
Let us set $R_i^0=N(X_i)$ and $R_i^1=\zet {\ka B_i}$, $i\in [d]$. Then $\Modgen{{X_i}}=R_i^0\oplus R_i^1$ and
$$
\Modgen X=\bigoplus_{\varepsilon\in\zet_2^d} R_1^{\varepsilon_1}\otimes\cdots\otimes R_d^{\varepsilon_d}
=\zet{\ka B}\oplus\bigoplus_{\varepsilon\neq (1,\ldots ,1)} R_1^{\varepsilon_1}\otimes\cdots\otimes R_d^{\varepsilon_d}.
$$
Now, it follows from Theorem \ref{rzut} that
$$
N(X)= \bigoplus_{\varepsilon\neq (1,\ldots ,1)}R_1^{\varepsilon_1}\otimes\cdots\otimes R_d^{\varepsilon_d}.
$$
This equation can be the launching point for an alternative approach to our theory.
It seems to be even simpler than ours but at the same time less natural. There is yet
another approach which we found at the beginning of our investigations. We are going
to publish it elsewhere.
\end{uwa}

\section{Words}
\label{words}
Let $S$ be a non-empty set called an \textit{alphabet}. The elements of $S$ will be called \textit{letters}. A permutation $s\mapsto s'$ of the alphabet $S$ such that $s''=(s')'=s$ and $s'\neq  s$ is said to be a \textit{complementation}. Each sequence of letters $s_1\cdots s_d$ is called a \textit{word} of length $d$. The set of all words of length $d$ is denoted by $S^d$. Two words $s_1\cdots s_d$ and $t_1\cdots t_d$ are \textit{dichotomous} if there is an $i\in[d]$
such that $s'_i=t_i$. 

In connection with applications to the cube tilings, it is suitable to consider $d$-boxes of arbitrary cardinality. From now $X$ is a $d$-box if $X=X_1\times\cdots\times X_d$ and $|X_i|>2$ for each $i\in [d]$. Definitions of a proper box, a suit, a polybox etc. remain unchanged. Suppose now that for each $i\in[d]$, we have  a mapping $f_i:S\to \pud {X_i}$ such that $f_i(s')=X_i\setminus f_i(s)$. Then we can define the mapping $f\colon S^d\to \pud  X$ by
$$
f(s_1\cdots s_d)=f_1(s_1)\times\cdots\times f_d(s_d).
$$ 
There is an obvious parallelism between $f$ and the mappings that preserve dichotomies. Therefore, we shall refer to $f$ as a mapping that \textit{preserves dichotomies} as well. If $W\subseteq S^d$, then $f(W)=\{f(w)\colon w\in W\}$ is said to be a \textit{realization} of the set of words $W$. The realization is said to be \textit{exact} if for each pair of words  
$v=s_1\cdots s_d$, $w=t_1\cdots t_d$ in $W$,  if  $s_i\not\in \{t_i, t'_i\}$, then $f_i(s_i)\not\in\{f_i(t_i), X_i\setminus f_i(t_i)\}$. 

If $W$ consists of pairwise dichotomous words, then we call it a \textit{(polybox) genome}. Each realization of a genome is a proper suit. Two genomes $V$ and $W$ are \textit{equivalent} if the realizations $f(W)$ and $f(V)$ are suits of the same polybox for each  mapping $f$ that preserves dichotomies. Now we collect several criteria of the equivalence. Two of them have been already discussed in previous sections, where they are expressed in terms of suits.

We begin with the criterion which is a consequence of Theorem \ref{rzut} and is described in detail in Remark \ref{kryt2}.
Let us expand the alphabet $S$  by adding an extra element $*$. This new set of symbols will be denoted by $*S$. Let $S_+$ and $S_-$ be arbitrary subsets of $S$ satisfying the following equations 
$$
S_-=\{s'\colon s\in S_+\}=S\setminus S_+.
$$
Let $*S_+=S_+\cup \{*\}$. Let us consider the ring $\zet[*S]$ over $\zet$ freely generated by $*S$. Let us identify each $s\in S_-$ with $*-s'$, and denote the ring which we obtain by this identification from $\zet[*S]$ by $R$. It is clear that each word $v\in S^d$ can be expanded in $R$ so that it is written in a unique way in the form of an element $v_+$ of the free module $\zet{*S_+}\subset R$. For  a finite $V\subset S^d$,  we define  $V\oplus \in \zet{*S_+}$ by the equation $V^\oplus =\sum_{v\in V} v_+$.                 

\begin{tw}
\label{equigenom1}
Let $V$, $W\subset S^d$ be two genomes. Then $V$ and $W$ are equivalent if and only if $V^\oplus =W^\oplus$.
\end{tw}

For each $u\in (*S)^d$, let the mapping $\varphi_u\colon  S^d\to \er$ be defined by the formula
$$
\varphi_u(w)=\prod_{i=1}^d ([s_i=t_i]-[s_i=t'_i]+[s_i=*]),
$$
where we let $u=s_1\cdots s_d$, $w=t_1\cdots t_d$.
We can define the index for genomes being the counterpart of the index for suits 
$$
\Index W u=\sum_{w\in W} \varphi_u(w).
$$

Theorem \ref{kryt1} leads to the  following result.
\begin{tw}
\label{equigenom2}
Two genomes $V$, $W\subset S^d$ are equivalent if and only if  for each $u\in (*S)^d$ 
$$
\Index V u =\Index W u.
$$
\end{tw}

Let us remark that in order to check the equivalence we do not consider all $u$, in fact it suffices to restrict ourselves to $u\in (*S_+)^d$ or even to those words $u=s_1\cdots s_d$, that for each $i\in[d]$ at least one of the letters $s_i$, $s_i'$  appears at $i$-th place of a certain word from $V\cup W$.     

Let $W\subset S^d$ be a genome and let $v\in S^d$ be a word. We say that $v$ is \textit{covered} by $W$, and write $v\sqsubseteq W$, if $f(v)\subseteq \bigcup f(W)$ for every mapping $f$ that preserves dichotomies. Let us define 
$g\colon S^d\times S^d\to \zet$ by the formula
$$
g(v,w)=\prod^d_{i=1}(2[s_i=t_i]+[s_i\not\in\{t_i, t'_i\}]).     
$$

Let 
$$
 E(S)=\{B\subset S\colon |\{s,s'\}\cap B|=1, \text{whenever $s\in S$}\}.
$$  
Similarly as in Section \ref{rikitiki}, for each letter $s\in S$, let us put $E s=\{B\in E(S)\colon s\in B\}$.
Lemma \ref{intermediolan} can be rephrased as follows

\begin{lemat}
\label{rikitiki2}
If $S$ is finite, $s_1,\ldots , s_n\in S$ and $s_i\not\in\{s_j, s'_j\}$ for every $i\neq j$, then
$$
|E s_1 \cap\dots\cap E s_n|=\frac{1}{2^{n}}|E(S)|.
$$
\end{lemat} 

For each word $u=s_1\cdots s_d\in S^d$, let us define $\breve{u}=Es_1\times\cdots\times  Es_d$.
It is clear that the mapping $S^d\ni u\mapsto \breve{u}$ preserves dichotomies.

\begin{tw} 
\label{Juventus}
Let $v\in S^d$ and let $W\subset S^d$ be a genome. Then $\sum_{w\in W} g(v,w)\le 2^d$. Moreover,
the following statements are equivalent:
\begin{itemize}
\item[{\em(1)}] $v \sqsubseteq W$, 
\item[{\em(2)}] $\breve{v}\subseteq \bigcup_{w\in W}\breve{w}$,
\item[{\em(3)}] $\sum_{w\in W} g(v,w)= 2^d$.

\end{itemize}
\end{tw}

\proof
Since the number of words involved is finite, the number of letters that constitute these words is finite. Therefore, we may restrict ourselves to a finite  subset of the alphabet $S$. (Such a restriction will influence the sets $\breve{u}$, $u\in W\cup \{v\}$, as they depend on the alphabet, but it does not influence the relation (2).) Observe that by Lemma \ref{rikitiki2} and the definition of $g$ one has
$$ 
|\breve{v}\cap\breve{w}|=\frac{g(v,w)|\breve{v}|}{ 2^d}.
$$    
Since $\breve{w}$, $w\in W$, form a suit, we get 
\begin{equation}
\label{brewa}
|\breve{v}|\ge |\breve{v}\cap\bigcup_{w\in  W}\breve{w}|=\sum_{w\in W}\frac{g(v,w)|\breve{v}|}{2^d}, 
\end{equation}
which immediately implies the first part of our theorem. If $v\sqsubseteq W$, then $\breve{v}\subseteq \bigcup_{w\in W}\breve{w}$ and (\ref{brewa}) becomes an equation. Therefore, it remains to show, that if $\sum_{w\in W}g(v,w)=2^d$,
then $v\sqsubseteq W$. Equivalently, if  $\breve{v}\subseteq \bigcup_{w\in W}\breve{w}$, then $v\sqsubseteq W$. Suppose that $v\not\sqsubseteq W$. Then there is a $d$-box $X$ and  a mapping $f\colon S^d\to \pud X$ that preserves dichotomies such that $f(v)\not\subseteq \bigcup f(W)$. 
 
\medskip
\noindent\textbf{Claim}\, If $Y$ is a 1-box and $h\colon S\to \pud Y$ preserves complementarity, that is, $h(s')=Y\setminus h(s)$ for  each $s\in S$, then for every $y\in Y$ there is a $C\in E(S)$ such that 
$$
[y\in h(s)]=[s\in C].
$$

\medskip
It suffices to put $C=\{s\colon [y\in h(s)]=1\}$. 

Suppose that $x\in f(v)\setminus \bigcup f(W)$. By Claim  and the definition of $f$, there is a $b\in E(S)^d$, $b=(B_1,\ldots, B_d)$, such that for each $u\in S^d$, $u=s_1\cdots s_d$,
$$
[x\in f(u)]=\prod_{i=1}^d[x_i\in f_i(s_i)]=\prod_{i=1}^d[u_i\in B_i]=[b\in \breve{u}]. 
$$
Therefore, $b\in \breve{v}\setminus \bigcup_{w\in W}\breve{w}$ and consequently $\breve{v}\not\subseteq \bigcup_{w\in W}\breve{w}$.
\hfill $\square$    

\medskip
Let $V$, $W\subset S^d$ be two genomes. We say that $W$ \textit{covers} $V$, and write $V\sqsubseteq W$, if each $v\in V$ is covered by $W$.

\begin{pr}
\label{rownopod}
Two genomes $V$ and $W$ contained in $S^d$ are equivalent if and only if $W$ covers $V$ and $|V|=|W|$.
\end{pr}

\proof The implication `$\Rightarrow$' is obvious.

We may assume that $S$ is finite. Then $|\breve{u}|=2^{(|S|-2)d}$ for each $u\in S^d$. Since $|V|= |W|$, we get
$$ 
|\bigcup_{v\in V}\breve{v}|=2^{(|S|-2)d}|V|=|\bigcup_{w\in W}\breve{w}|.
$$
This equation together with the assumption $V\sqsubseteq W$ give us $\bigcup_{v\in V}\breve{v}= \bigcup_{w\in W}\breve{w}$.
Therefore, for each $w\in W$ we have $\breve{w}\subseteq \bigcup_{v\in V}\breve{v}$ which, by Theorem \ref{Juventus}, implies $w\sqsubseteq V$.\hfill $\square$   

\begin{uwa} 
By analogy, we can define binary codes for $S^d$. One can deduce from Proposition \ref{rownopod} that $V$ and $W$ are equivalent if and only if $\beta (V)=\beta (W)$  for each binary code $\beta$.  
\end{uwa}

\medskip
For $w\in S^d$, $w=s_1\cdots s_d$, and $\varepsilon \in \zet_2^d$, define 
$w^\varepsilon=s_1^{\varepsilon_1}\cdots s_d^{\varepsilon_d}$ in the same way as it has been done for boxes. Denote by $C_w$ the set $\{w^\varepsilon\colon \varepsilon\in \zet_2^d\}$.  
An inspection of the argument used in the proof of Theorem \ref{sztywnypal} suggests the following 
\begin{tw}
\label{niejednoznacznosc}
Let $v\in S^d$, and $W\subset S^d$ be a genome. If $v\sqsubseteq W$ and $v\not\in W$, then there is a 
$w\in W$ such that 
$$
|\Index W w|<|C_w\cap W|.
$$
\end{tw}

Let $C=C_w$ for some $w\in S^d$. As in the case of simple suits, we say that $u$, $v\in C$ carry the same sign if $(-1)^{|\varepsilon|}=1$ for the only $\varepsilon$ such that $u^\varepsilon=v$. The orientation $C^-$, $C^+$ of $C$ is defined accordingly (cf. Section \ref{indices}). Suppose that $W$ is a genome and suppose that for each $w\in W$ we have chosen an orientation $C_w^+$, $C_w^-$ of $C_w$. These orientations define the following decomposition of $W$: $W^+=W\cap \bigcup_{w\in W}C^+_w$, $W^-=W\cap \bigcup_{w\in W}C^-_w$. We call it \textit{induced}.

\begin{tw}[rigidity of genomes] 
\label{rikitikigenom}
If $W$, $V\subset S^d$ are equivalent genomes, $W^+$, $W^-$ is an induced decomposition of $W$ and  $W^+\subseteq V$, then $W=V$. 
\end{tw}

\proof Let $v\in V\setminus W^+$. Then $v\sqsubseteq W^-$. If $v\not\in W^-$, then by the preceding theorem there is a $w\in W^-$ such that  
$
|\Index {W^-} w|<|C_w\cap W^-|,
$  
which contradicts the definition of  $W^-$.
\hfill$\square$ 
\section{Rigidity of cube tilings}
\begin{pr}
\label{prpr}
Let $W\subset S^d$ be a genome. Let $X$ be a $d$-box whose cardinality can be infinite.   
If $f\colon S^d\to \pud X$ preserves dichotomies, then the realization $f(W)$ is a proper suit for $X$ if and only if $|W|=2^d$.
\end{pr}

\proof For each $X_i$ one can find a finite subset $Y_i$ such that $Y$ is a $d$-box and $f(w)\cap Y\in \pud Y$, whenever $w\in W$. If $f(W)$ is a proper suit for $X$, then $\ka F =\{f(w)\cap Y\colon w\in W\}$ is a proper suit  for $Y$. 
Therefore, $\ka F$ has to contain $2^d$ elements (see Section \ref{minimal}, also \cite[Theorem 2]{GKP}). On the other hand, $|\ka F|=|W|.$   
 
Suppose now that there are $W$ and $f$ such that $\bigcup f(W)\neq X$. Let $x\in X\setminus \bigcup f(W)$. Then there is a finite $d$-box $Y$ such that $x\in Y$ and $\ka F =\{f(w)\cap Y\colon w\in W\}\subset \pud Y$ is a proper suit. Since $|W|=2^d$, then, by Theorem \ref{krakow}, $\ka F$ is a partition of $Y$. Thus $x\in f(w)$ for some $w$, which is a contradiction.\hfill$\square$     

\medskip
Our next lemma follows immediately from the preceding proposition and the definition of the cover relation $\sqsubseteq$.
\begin{lemat}
\label{stoliczek}
If $W\subseteq S^d$ is a genome that consists of $2^d$ elements, then  $v\sqsubseteq W$ for each $v\in S^d$.
\end{lemat} 
\begin{tw}
\label{laga}
Let $W\subset S^d$ be a genome that consists of $2^d$ elements. Let $W^+$ and $W^-$ be an induced decomposition of $W$. If $\{v\}\cup W^+$ is a genome, then $v\in W^-$. 
\end{tw}

\proof By the above lemma and the fact that $\{v\}\cup W^+$ is   a genome, we deduce that $v\sqsubseteq W^-$. The conclusion follows now immediately from the definition of $W^-$ and Theorem \ref{niejednoznacznosc}. \hfill$\square$  
 
\begin{uwa}
This result verifies the \textit{rigidity conjecture for 2-extremal cube tilings} of Lagarias and Shor \cite{LS2} 
which states that 
$W^+$ determines $W^-$ uniquely if $|C_w\cap W|=2$ for each $w\in W$. 
\end{uwa}


\begin{pr}
\label{minimini}
If $X$ is a $d$-box of arbitrary cardinality and $\ka F$ is a partition of $X$ into proper boxes such that $|\ka F|=2^d$, then $\ka F$ is a suit. 
\end{pr}

\proof Fix two elements $A$, $B\in \ka F$. Pick a finite box $Y\subseteq X$ so that $\ka F_Y=\{ D\cap Y\colon D\in \ka F\}\subset \pud Y$. By definition,  $\ka F_Y$ is a partition of $Y$ which consists of $2^d$ elements. By (\ref{vol0}) and Theorem \ref{vol00}, $\ka F_Y$ is a minimal partition of $Y$. It follows now from  Theorem \ref{minor} that it is a suit. Hence there is an $i\in [d]$ such that $A_i\cap Y_i=Y_i\setminus B_i$. It is clear that if $A$ and $B$ were not dichotomous, then for properly chosen $Y$ their intersections with $Y$ would not be dichotomous in $Y$.
\hfill$\square$

\medskip
Let $T$ be a subset of $\mathbb{R}^d$. The family $[0,1)^d+T=\{I_t\colon t\in T\}$, where  $I_t=[0,1)^d+t$, is  a \textit{cube tiling} of $\mathbb{R}^d$ if $\mathbb{R}^d=\bigcup_{t\in T}I_t$ and $I_u\cap I_v=\emptyset$ for distinct $u,v \in T$. This cube tiling is said to be $2\mathbb{Z}^d$-\textit{periodic} if $T$ is $2\mathbb{Z}^d$-periodic. 

Consider the alphabet $S=(-1,1]$ with complementation $s\mapsto s'$ defined by the equation $|s-s'|=1$. We write the elements of $S^d$ as vectors: $t=(t_1,\ldots, t_d)$. Let $J_t=I_t\cap [0,1]^d$. Define 
$$
W=\{t\in T\colon J_t\neq\emptyset\}\subset S^d.
$$     
It is easily seen that each box from $\ka G=\{J_t\colon t\in W\}$ contains exactly one vertex of the cube $[0,1]^d$. Therefore $|\ka G|=2^d$. Theorem \ref{minimini} implies now that $\ka G$ is a suit. The latter fact leads to the conclusion that $W$ is a genome. Let $W^+$, $W^-$ be an induced decomposition of $W$. Define $T^+=W^++2\zet^d$ and 
$T^-=W^-+2\zet^d$.  
\begin{tw}[chess-board decomposition]
\label{krewmleko}
Let $[0,1)^d +T$ be a cube tiling of $\er^d$, $T^+$, $T^-$ be as defined  above and $z\in \er^d$. If $I_z$ is disjoint   
with $[0,1)^d+T^+$, then $z\in T^-$. 
\end{tw}

\proof
Let $K_z=[0,1]^d+z$. Define $U\subset (-1,1]^d$ so that $u\in U$ if and only if $I_{u+z}$ intersects $K_z$. As in the case $W$, the set $U$ is  a genome and $|U|=2^d$. Let $U^+=U\cap (T^+-z)$, $U^-=U\cap (T^--z)$. Clearly, $U^+$, $U^-$ is an induced decomposition of $U$. Observe that since $I_z$ is disjoint with  $[0,1)^d+T^+$, $0$ is dichotomous to each element of $U^+$. Therefore, by Theorem \ref{laga}, we have $0\sqsubseteq U^-$, which in turn implies $z\in T^-$.
\hfill$\square$

\end{document}